\documentclass{article}[10]

\usepackage{latexsym}
\usepackage{amssymb} 
\usepackage{epsfig}
\usepackage{amsmath}
\usepackage[mathscr]{eucal}

\newtheorem{Lemma}{Lemma}
\newtheorem{Proposition}[Lemma]{Proposition}
\newtheorem{Theorem}[Lemma]{Theorem}
\newtheorem{Example}{Example}

\newtheorem{Definition}{Definition}

\newcommand{\cd}{\ \stackrel{d}{\longrightarrow} \ }
\newcommand{\ed}{\ \stackrel{d}{=} \ }

\newcommand{\PP}{\mbox{${\mathcal P}$}}

\newcommand{\VV}{\mbox{${\mathcal V}$}}
\newcommand{\GG}{\mbox{${\mathcal G}$}}
\newcommand{\FF}{\mbox{${\mathcal F}$}}

\newcommand{\HH}{\mbox{${\mathcal H}$}}
\newcommand{\AAA}{\mbox{${\mathcal A}$}}

\newcommand{\II}{\mbox{${\mathcal I}$}}

\newcommand{\eps}{\varepsilon}
\newcommand{\bE}{{\bf E}}
\newcommand{\bP}{{\bf P}}

\newcommand{\bi}{{\bf i}}

\newcommand{\TTT}{\mbox{${\mathscr{T}}$}}
\newcommand{\EEE}{\mbox{${\mathscr{E}}$}}
\newcommand{\VVV}{\mbox{${\mathscr{V}}$}}

\newcommand{\VVVV}{\mbox{$\widetilde{\VVV}$}}

\newcommand{\var}{{\bf Var}}

\newcommand{\sfrac}[2]{{\textstyle\frac{#1}{#2}}}

\newcommand{\qed}{\hfill{\ \ \rule{2mm}{2mm}}}
\newcommand{\proof}{\noindent\emph{Proof :}\ \ \ \ }

\newcommand{\mudiag}{\mu^\nearrow}

\newcommand{\lleq}{\leq^*}

\newcommand{\muprod}{\mu \otimes \mu}
\newcommand{\Tprod}{T \otimes T}
\newcommand{\norm}{\parallel}

\newcommand{\Rbold}{{\mathbb{R}}}
\newcommand{\Tbold}{{\mathbb{T}}}
\newcommand{\Zbold}{{\mathbb{Z}}}
\newcommand{\Nbold}{{\mathbb{N}}}

\title{{\bf A Necessary and Sufficient Condition for the Tail-Triviality
      of a Recursive Tree Process} 
\footnote{\emph{Abbreviated title} : Condition for tail-triviality of a RTP}} 

\author{{\bf Antar Bandyopadhyay} 
\footnote{E-Mail : {\tt antar@math.chalmers.se}. 
This work started while the author was
a gradute student at the Department of Statistics, University of
California, Berkeley.}
\vspace{0.1in} \\ 
        Department of Mathematics \\
        Chalmers University of Technology \\
        SE-412 96 G\"{o}teborg, SWEDEN}

\begin{document}


\maketitle

\begin{abstract}

Given a \emph{recursive distributional equation} (RDE) and a solution $\mu$
of it, we consider the tree indexed invariant process called the
\emph{recursive tree process} (RTP) with marginal $\mu$.
We introduce a new type of bivariate uniqueness property which is
different from the one defined by Aldous and Bandyopadhyay \cite{AlBa05}, 
and we prove that this property is equivalent to tail-triviality for the
RTP, thus obtaining a necessary and sufficient condition to determine
tail-triviality for a RTP in general. As an application we consider
Aldous' construction of the frozen percolation process on a infinite
regular tree \cite{Al00} and show that the associated RTP has a trivial tail.

\end{abstract}

\vspace{0.25in}


\ \\


\newpage

\section{Introduction, Background and Motivation}
\label{Sec:Intro}
Fixed-point equations or distributional identities have appeared in the
probability literature for quite a long time in a variety of settings. 
The recent survey of Aldous and Bandyopadhyay \cite{AlBa05}
provides a general framework to study certain type of distributional 
equations. 

Given a space $S$ write $\PP\left(S\right)$ for the set of all probabilities
on $S$. A \emph{recursive distributional equation} (RDE) \cite{AlBa05} is a 
fixed-point equation on $\PP\left(S\right)$ defined as
\begin{equation}
X \ed g \left( \xi ; \left( X_j : 1 \leq j \lleq N \right) \right)
\,\,\,\, \mbox{on} \,\,\, S, 
\label{Equ:RDE}
\end{equation}
where it is assumed that
$\left(X_j\right)_{j \geq 1}$ are i.i.d. $S$-valued random variables
with same distribution as
$X$, and are independent of the pair $\left( \xi, N \right)$. Here
$N$ is a non-negative integer valued random variable, 
which may take the value $\infty$, and $g$ is a given $S$-valued function. 
(In the above equation by ``$\lleq N$'' we mean the left hand side is 
``$ \leq N$'' if $N < \infty$, and ``$< N$'' otherwise). 
In (\ref{Equ:RDE}) the distribution of $X$ is \emph{unknown} while
the distribution of the pair $\left( \xi, N \right)$, and the function
$g$ are the \emph{known} quantities. Perhaps a more conventional 
(analytic) way of writing the equation (\ref{Equ:RDE}) would be
\begin{equation}
\mu = T\left(\mu\right) \, , 
\label{Equ:RDE-T}
\end{equation}
where $T : \PP \rightarrow \PP\left(S\right)$ is a function defined on
$\PP \subseteq \PP\left(S\right)$ such that $T\left(\mu\right)$ is
the distribution of the right-hand side of the equation (\ref{Equ:RDE}), 
when $\left(X_j\right)_{j \geq 1}$ are i.i.d. $\mu \in \PP$. 

As outlined in \cite{AlBa05} 
in many applications RDEs play a very crucial 
role. Examples include study of Galton-Watson branching processes and
related random trees, probabilistic analysis of algorithms with
suitable recursive structure \cite{RosRu01, Ros92}, 
statistical physics models on trees \cite{AlSt04, Al00, Ga04, Ba05a},
and statistical physics and algorithmic questions in the mean-field model
of distance \cite{Al92a, Al01, AlSt04}.
In many of these applications, particularly in the last 
two types mentioned above, often one needs to construct 
a particular tree indexed \emph{stationary}
process related to a given RDE, which is called a
\emph{recursive tree process} (RTP) \cite{AlBa05}.

\subsection{Recursive Tree Process}
\label{SubSec:RTP}  
More precisely, suppose the RDE (\ref{Equ:RDE}) has a 
solution, say $\mu$. Then as shown in \cite{AlBa05}, using the
consistency theorem of Kolmogorov \cite{Bill95},
one can construct a process, say 
$\left(X_{\bi}\right)_{\bi \in \VV}$, indexed by 
$\VV := \left\{ \emptyset \right\} \cup_{d \geq 1} \Nbold^d$, such that 
\begin{equation}
\begin{array}{cl}
\mbox{(i)} & X_{\bi} \sim \mu \,\,\,\, \forall \,\,\, \bi \in \VV, \\
\mbox{(ii)} & 
\mbox{For each\ \ } d \geq 0, \left(X_{\bi}\right)_{\vert \bi \vert = d}
\mbox{\ \ are independent}, \\
\mbox{(iii)} & 
X_{\bi} = g \left( \xi_{\bi} ; \left( X_{\bi j} : 
1 \leq j \lleq N_{\bi} \right) \right) \,\,\,\, \forall \,\,\, \bi \in \VV, \\
\mbox{(iv)} &
X_{\bi} \mbox{\ \ is independent of\ \ } 
\left\{ \left(\xi_{\bi'}, N_{\bi'} \right) \,\Big\vert\, 
\vert \bi' \vert < \vert \bi \vert \, \right\}
\,\,\,\, \forall \,\,\, \bi \in \VV,
\end{array}
\label{Equ:RTP}
\end{equation}
where $\left(\xi_{\bi}, N_{\bi}\right)_{\bi \in \VV}$ are taken to be i.i.d. 
copies of the pair $\left( \xi, N \right)$, and by 
$\vert \cdot \vert$ we mean the length of a finite word. 
The process $\left(X_{\bi}\right)_{\bi \in \VV}$ is called an invariant 
\emph{recursive tree process} (RTP) with marginal $\mu$. The i.i.d.
random variables $\left(\xi_{\bi}, N_{\bi} \right)_{\bi \in \VV}$ are 
called the \emph{innovation process}. In some sense an invariant 
RTP with marginal $\mu$, 
is an almost sure representation of a solution $\mu$ of the RDE
(\ref{Equ:RDE}). 
Here we note that there is a natural tree structure on $\VV$. 
Taking $\VV$ as the vertex set, we join two words $\bi, \bi' \in \VV$
by an edge, if and only if, $\bi' = \bi j$
or $\bi = \bi' j$, for some $j \in \Nbold$. We will denote this
tree by $\Tbold_{\infty}$. 
The empty-word $\emptyset$ will
be taken as the root of the tree $\Tbold_{\infty}$, and we will write
$\emptyset j = j$ for $j \in \Nbold$. 

In the applications mentioned above
the variables $\left(X_{\bi}\right)_{\bi \in \VV}$ 
of a RTP are often used as auxiliary variables to define or to construct some
useful random structures. To be more precise in \cite{Al01} they were used
to obtain ``almost optimal matching'', while in \cite{Al00} they were
used to define the percolation clusters. 
In such applications typically the innovation process defines the
``internal'' variables while the RTP is constructed ``externally'' using
the consistency theorem. It is then natural
to ask whether the RTP is measurable only with respect to 
the i.i.d. innovation process $\left(\xi_{\bi}, N_{\bi}\right)$. 
\begin{Definition}
\label{Defi:Endogeny}
An invariant RTP with marginal $\mu$ is called \emph{endogenous}, if the
root variable $X_{\emptyset}$ is measurable with respect to the
$\sigma$-algebra 
\[ \GG := \sigma \left( \left\{ \left(\xi_{\bi}, N_{\bi} \right)
\,\Big\vert\, \bi \in \VV \,\right\} \right). \]
\end{Definition}
This notion of endogeny has been the main topic of discussion in
\cite{AlBa05}. The authors provide a necessary and sufficient
condition for endogeny in the general setup \cite[Theorem 11]{AlBa05}.
A non-trivial application of this result is given in \cite{Ba02},
where it is proved that
the invariant RTP associated with the \emph{logistic} RDE, which appears
in the study of the mean-field random assignment problem \cite{Al01} is
endogenous. Another interesting example arise in the construction of the
\emph{frozen percolation} on an infinite $3$-regular tree 
by Aldous \cite{Al00}, where a particular RTP has been used to carry on
the construction. This example is one of our main motivations, so
we discuss this example in more detail in Section \ref{SubSec:Appl-FPP}.

As discussed in \cite{AlBa05} in some sense, the concept of endogeny 
tries to capture the idea of having ``no influence of the boundary at 
infinity'' on the root. In this direction a closely 
related concept would be the tail-triviality of a RTP. To give a formal 
definition of the tail of a RTP, let
$\left(X_{\bi}\right)_{\bi \in \VV}$ be an invariant RTP with 
marginal $\mu$, where $\mu$ is a solution of the RDE (\ref{Equ:RDE}). 
The tail $\sigma$-algebra of $\left( X_{\bi} \right)_{\bi \in \VV}$ is
defined as 
\begin{equation} 
\HH = \mathop{\cap}\limits_{n \geq 0} \HH_n,
\label{Equ:RTP-Tail}
\end{equation}
where
\begin{equation}
\HH_n := \sigma \left( \left\{ X_{\bi} \,\Big\vert\, \vert \bi \vert 
\geq n \,\right\} \right).
\label{Equ:RTP-Tail-n}
\end{equation}
Naturally, we will say an invariant RTP has trivial tail if the tail
$\sigma$-algebra $\HH$ is trivial. 
Because the innovation process $\left(\xi_{\bi}, N_{\bi}\right)_{\bi \in \VV}$
is i.i.d., so it is natural to expect that if a RTP is endogenous, then
it has a trivial tail.
\begin{Proposition} 
\label{Prop:Endo}
Suppose $\mu$ is a solution of the RDE (\ref{Equ:RDE}) and 
$\left( X_{\bi} \right)_{\bi \in \VV}$ be an invariant RTP with marginal
$\mu$. Then the tail of $\left(X_{\bi}\right)_{\bi \in \VV}$ is trivial 
if it is endogenous. 
\end{Proposition}
Thus one way to conclude that a RTP is not endogenous will be to show that
it has a non-trivial tail. The following easy example shows that the
converse may not hold. 
\begin{Example}
\label{Ex:MC}
Take $S := \{ 0, 1 \}$. Let $0 < q < 1$ and 
$\xi \sim \mbox{Bernoulli}\left(q\right)$. Consider the RDE
\begin{equation}
X \ed \xi + X_1 \,\,\,\, \left( \, \mbox{\emph{mod}\ \ } 2 \, \right),
\label{Equ:RDE-Ex-MC}
\end{equation}
where $X_1$ has same distribution as $X$, and is independent of $\xi$. 
\end{Example}
If $T$ is the associated operator defined by the right-hand side of the
equation (\ref{Equ:RDE-Ex-MC}), then it is easy to see that $T$ maps a
$\mbox{Bernoulli}\left(p\right)$ distribution to a 
$\mbox{Bernoulli}\left(p'\right)$ distribution where
\[ p' = p \left( 1 - q \right) + q \left( 1 - p \right). \]
Thus the unique solution of the RDE (\ref{Equ:RDE-Ex-MC}) is 
$\mbox{Bernoulli}\left(\sfrac{1}{2}\right)$. 

In this example because there is no branching 
$\left( N \equiv 1 \right)$, so 
the invariant RTP with marginal $\mbox{Bernoulli}\left(\sfrac{1}{2}\right)$
can be indexed by the non-negative integers, we denote it by 
$\left(X_i\right)_{i \geq 0}$, where $X_0$ is the root variable and it
satisfy
\[ 
X_i = \xi_i + X_{i+1} \,\,\,\, \mbox{a.s.} \,\,\, \forall \,\,\, i \geq 0,
\]
where $\left( \xi_i \right)_{i \geq 0}$ are i.i.d. 
$\mbox{Bernoulli}\left(q\right)$. It is then easy to see that we must have 
\[
X_{i+1} \mbox{\ and\ } \left( \xi_0, \xi_1, \ldots, \xi_i \right)
\,\,\,\, \mbox{are independent, for all\ \ } i \geq 0.
\]
Therefor $X_0$ is independent of the innovation process
$\left( \xi_i \right)_{i \geq 0}$, thus the RTP is not endogenous. 
The following proposition whose proof we defer till 
Section \ref{Sec:Root-Tail}, states that the RTP
$\left( X_i \right)_{i \geq 0}$ has trivial tail. This gives an example of 
an invariant RTP which is not endogenous but has trivial tail. 
\begin{Proposition}
\label{Prop:Ex-MC}
The invariant RTP with marginal $\mbox{Bernoulli}\left(\sfrac{1}{2}\right)$
associated with the RDE (\ref{Equ:RDE-Ex-MC}) has trivial tail. 
\end{Proposition}

So proving tail-triviality of a RTP is weaker than proving endogeny, but 
in some cases it might help to prove non-endogeny by showing that the tail
is not trivial. Also in general,
studying the tail of a stochastic process is 
mathematically interesting. 

In this article we provide a necessary and sufficient condition 
to determine the tail-triviality for an invariant RTP. This condition is in
the same spirit of the equivalence theorem 
of Aldous and Bandyopadhyay \cite[Theorem 11]{AlBa05}. 
But before we state our main result we first introduce a new type of
\emph{bivariate uniqueness} property, which is different than the
one introduced in \cite{AlBa05}, we will call it the
\emph{bivariate uniqueness property of the second kind}. 

\subsection{Bivariate Uniqueness Property of the Second Kind}
\label{Subsec:Bi-uni-2nd}
Consider a general RDE given by (\ref{Equ:RDE}) and let 
$T \colon \PP \rightarrow \PP\left(S\right)$ be the 
induced operator. We will consider a bivariate version of it. Write
$\PP^{(2)}$ for the space of probability measures
on $S^2 = S \times S$, with marginals in $\PP$. 
We can now define 
a map $\Tprod: \PP^{(2)} \rightarrow \PP\left(S^2\right)$ as follows
\begin{Definition}
\label{Def:Bi-T-2nd}
For a probability $\mu^{(2)} \in \PP^{(2)}$, 
$\left(\Tprod\right)\left(\mu^{(2)}\right)$ is the joint distribution of 
\[
\left( \begin{array}{c}
       g\left( \xi, X_j^{(1)}, 1 \leq j \lleq N \right) \\
       g\left( \eta, X_j^{(2)}, 1 \leq j \lleq M \right)
       \end{array} \right)
\]
where we assume 
\begin{enumerate}
\item $\left(X_j^{(1)}, X_j^{(2)}\right)_{j \geq 1}$ are
      independent with joint distribution $\mu^{(2)}$ on $S^2$; 
\item $\left(\xi, N\right)$ and 
      $\left(\eta, M\right)$ are i.i.d;
\item the families of random variables in 1 and 2 are independent. 
\end{enumerate}
\end{Definition}
We note that here we use \emph{independent} copies of the innovation pair
in the two coordinates. We also note that this is preciously where 
this bivariate operator differs from the bivariate operator 
defined in \cite{AlBa05}, where the innovation pair was kept same
at each coordinate.

From the definition it follows that 
\begin{Lemma}
\label{Lem:Bi-Obvious-Sol}
\begin{itemize}
\item[(a)] If $\mu$ is a fixed point for $T$, then the associated
           \emph{product measure} $\muprod$ is a fixed point 
           for $\Tprod$.
\item[(b)] If $\mu^{(2)}$ is a fixed point for $\Tprod$, then each marginal
           distribution is a fixed point for $T$.
\end{itemize}
\end{Lemma}
So if $\mu$ is a fixed point for $T$ then $\muprod$ is a fixed
point for $\Tprod$ and there may or may not be other fixed points
of $\Tprod$ with marginal $\mu$.
\begin{Definition}
\label{Def:Bi-equiv}
An invariant RTP with marginal $\mu$ has the
\emph{bivariate uniqueness property of the second kind} 
if $\muprod$ is the unique 
fixed point of $\Tprod$ with marginal $\mu$.
\end{Definition}

\subsection{Main Result : An Equivalence Theorem}
\label{Subsec:Equiv-Thm-2}
Our main theorem is the following general result linking 
the tail triviality of an invariant RTP with the
bivariate uniqueness property of the second  kind. 
\begin{Theorem}
\label{Thm:Bi-equiv-2}
Suppose $S$ is a \emph{Polish} space. 
Consider an invariant RTP with marginal distribution $\mu$.
\begin{itemize}
\item[(a)] If the RTP has trivial tail then the bivariate
           uniqueness property of the second  kind holds.
\item[(b)] Suppose the bivariate uniqueness property of the second  kind holds.
           If also $\Tprod$ is continuous with respect to weak
           convergence on the set of bivariate distributions with
           marginals $\mu$, then the tail of the RTP is trivial. 
\item[(c)] Further, the RTP has a trivial tail if and only if 
           \[ \left(\Tprod\right)^n \left( \mudiag \right) \cd \muprod, \]
           where $\mudiag$ is the diagonal measure with marginal $\mu$, that 
           is, if $\left( X, Y \right) \sim \mudiag$, then
           $\bP\left(X=Y\right) = 1$ and $X, Y \sim \mu$. 
\end{itemize}
\end{Theorem}

\subsection{Heuristic Behind the Equivalence Theorem}
\label{SubSec:Heuristic}
\begin{figure}[h]
\begin{center}
\epsfig{file=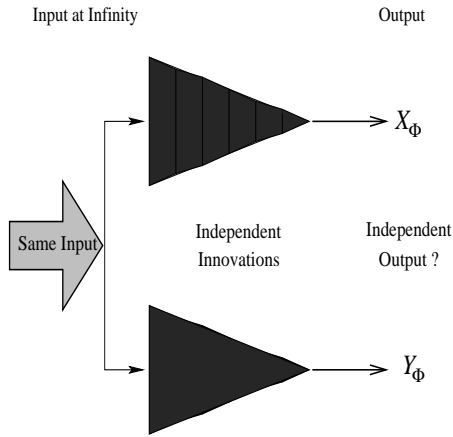, width=0.6\textwidth, height=0.6\textwidth}
\caption{Intuitive picture for the bivariate uniqueness of the second kind} 
\label{Fig:Theory-bi-uni-2}
\end{center}
\end{figure}
\noindent
Suppose $\mu$ is a solution of the RDE (\ref{Equ:RDE}) and let 
$\left(X_{\bi}\right)_{\bi \in \VV}$ be an invariant RTP with marginal 
$\mu$. Let $\left(\xi_{\bi}, N_{\bi}\right)_{\bi \in \VV}$ be the
i.i.d. innovation process, and 
$\GG_n := \sigma \left( \left\{ \left(\xi_{\bi}, N_{\bi} \right)
\,\Big\vert\, \vert \bi \vert \leq n \, \right\} \right)$ be the
$\sigma$-algebra for the innovations in first $n$-generations of the
tree $\Tbold_{\infty}$. From the construction (\ref{Equ:RTP}) of the
RTP we note that for any $n \geq 0$ the root variable $X_{\emptyset}$ is
measurable with respect to the the $\sigma$-algebra
$\sigma \left( \GG_n \cup \HH_{n+1} \right)$, where $\HH_{n+1}$ is
as defined in (\ref{Equ:RTP-Tail-n}). 
So heuristically to check whether the tail of the RTP 
$\HH = \cap_{n \geq 0} \HH_n$ contains any non-trivial information, we
may want to do the following :
\begin{quote}
Start with ``same input at infinity.'' Take two independent but 
identical copies of the innovation process and run 
through the recursions in (\ref{Equ:RTP}). Finally obtaining two 
copies of the RTP, say $\left(X_{\bi}\right)_{\bi \in \VV}$ and
$\left(Y_{\bi}\right)_{\bi \in \VV}$,
with same marginal $\mu$. Check if the root
variables $X_{\emptyset}$ and $Y_{\emptyset}$, are independent or not. 
\end{quote}
Figure \ref{Fig:Theory-bi-uni-2} gives this intuitive picture. 
The part (c) of the Theorem \ref{Thm:Bi-equiv-2} makes this 
process rigorous. Moreover we notice from definition the bivariate
process $\left(X_{\bi}, Y_{\bi}\right)_{\bi \in \VV}$ is a RTP 
associated with the operator $T \otimes T$. This leads to the notion of
bivariate uniqueness property of the second kind. We would like to note 
that the
proof of the Theorem \ref{Thm:Bi-equiv-2} is nothing but to make this
heuristic rigorous. 

\subsection{Application to Frozen Percolation}
\label{SubSec:Appl-FPP}
As mentioned earlier, one of our main motivating example arise in the
context of frozen percolation process on an infinite regular tree. For
sake of completeness we here provide a very brief background on
frozen percolation process, readers are advised to look at 
\cite{Al00, AlBa05} for more details. 

Frozen percolation process was first studied by Aldous \cite{Al00} where he 
constructed the process on a infinite $3$-regular tree. Let 
$\TTT_3 = \left( \VVV, \EEE \right)$ denote the infinite $3$-regular tree. 
Let $\left( U_e \right)_{e \in \EEE}$ be i.i.d $\mbox{Uniform}[0,1]$ 
edge weights. Consider a collection of random subsets 
$\AAA_t \subseteq \EEE$ for $0 \leq t \leq 1$, whose evolution is 
described informally by :
\begin{quote}
$\AAA_0$ is empty; for each $e \in \EEE$, at time 
$t = U_e$ set $\AAA_t = \AAA_{t-} \cup \{ e \}$ if
each end-vertex of $e$ is in a finite cluster of $\AAA_{t-}$;
otherwise set $\AAA_t = \AAA_{t-}$. \hfill $(*)$ 
\end{quote}
(A {\em cluster} is formally a connected component of edges, but we also
consider it as the induced set of vertices). Qualitatively, in the process 
$\left( \AAA_t \right)$ the clusters may grow to infinite size
but, at the instant of becoming infinite they are ``frozen'', in the
sense that no extra edge may be connected to an infinite
cluster. The final set $\AAA_1$ will be a forest on 
$\TTT_3$ with both infinite and finite clusters, such that no two finite
clusters are separated by a single edge. 
Aldous \cite{Al00} defines this process $\left( \AAA_t \right)$ as the
\emph{frozen percolation process}. 

Although this process is intuitively quite natural,
rigorously speaking it is
not clear that it exists or that $(*)$ does specify a unique process. In fact
Itai Benjamini and Oded Schramm have an argument that such
a process does not exist on the $\Zbold^2$-lattice
(see the remarks in Section 5.1 of \cite{Al00}).
But for the infinite $3$-regular tree,
\cite{Al00} gives a rigorous construction of an automorphism invariant
process satisfying $(*)$. This construction uses the following RDE
\begin{equation}
Y \ed \Phi \left( Y_1 \wedge Y_2 ; U \right) 
\,\,\,\, \mbox{on} \,\,\, I := \left[\sfrac{1}{2}, 1\right] \cup \{\infty\},
\label{Equ:FPP-RDE}
\end{equation}
where $\left(Y_1, Y_2\right)$ are i.i.d with same distribution as $Y$, and are
independent of $U \sim \mbox{Uniform}[0,1]$, and 
$\Phi : I \times [0,1] \rightarrow I$  is a function defined as
\begin{equation}
\Phi(x ; u) = \left\{ \begin{array}{ll}
                                  x & \mbox{if\ \ } x > u \\
                                  \infty & \mbox{otherwise} 
                     \end{array} \right. .
\label{Equ:FPP-Phi-defi}
\end{equation}  
We will call (\ref{Equ:FPP-RDE}) the \emph{frozen percolation RDE}.

It turns out \cite{Al00} that the RDE (\ref{Equ:FPP-RDE}) has many 
solutions. In particular, solutions having no atom in
$\left[\sfrac{1}{2}, 1\right]$ are given by
\begin{equation}  
\nu_a(dx) = \sfrac{dx}{2 x^2} \,\,\,\, \sfrac{1}{2} < x < a; \,\,\,\,
\nu_a\left( \left\{\infty\right\} \right) = \sfrac{1}{2 a}, 
\label{Equ:FPP-Sols} 
\end{equation} 
where $a \in \left[\sfrac{1}{2}, 1\right]$, thus $\nu = \nu_1$
is the unique solution with support $I$. 

Notice that for the RDE (\ref{Equ:FPP-RDE})
$N \equiv 2$, so a RTP with marginal $\nu$ essentially lives on a
rooted binary tree, we will denote the vertex set in this case by 
$\VVVV$. Let $\left(Y_{\bi}\right)_{\bi \in \VVVV}$ be an invariant 
RTP with marginal $\nu$. Aldous' construction of the
frozen percolation process \cite{Al00} uses these \emph{externally} 
defined random variables $\left(Y_{\bi}\right)_{\bi \in \VVVV}$. 
We refer the readers to look at \cite{Al00} for the technical
details of this construction. 
Here we only mention briefly what is the significance of the RTP 
$\left(Y_{\bi}\right)_{\bi \in \VVVV}$. 
Let $e=(u,v)$ be an edge of the infinite regular binary tree $\TTT_3$  and let 
$\overrightarrow{e} = \overrightarrow{(u,v)}$ be a direction of it which
is from the vertex $u$ to vertex $v$. Naturally the directed edge
$\overrightarrow{e}$ has two directed edges coming out of it, which can be
considered as two children of it. Continuing in similar manner we notice 
that each directed edge 
$\overrightarrow{e}$ represent a rooted infinite binary tree, which is 
isomorphic to $\VVVV$, and the weights are defined appropriately using the
i.i.d Uniform edge weights $\left( U_e \right)$. 
If the frozen percolation process exists, then the time for the edge
$e$ to join to infinite along the subtree defined by $\overrightarrow{e}$
is given by the variable $Y_{\emptyset}$. More preciously, such time
should satisfy the distributional recursion (\ref{Equ:FPP-RDE}). 
However to prove the existence of the 
process such times are then \emph{externally} constructed using 
the RTP construction. 
Naturally it make sense to ask whether these variables can only 
be defined using the i.i.d $\mbox{Uniform}[0,1]$ edge weights
(see Remark 5.7 in \cite{Al00}), which is same as asking whether the RTP is
endogenous. 

\begin{Theorem}
\label{Thm:FPP-Tail-Trivial}
Any invariant \emph{recursive tree process} associated with the RDE
(\ref{Equ:FPP-RDE}) with marginal $\nu$ has trivial tail. 
\end{Theorem}
This result does not resolve the question of endogeny, but it proves that
the version of the frozen percolation process constructed by Aldous
\cite{Al00} on an infinite $3$-regular tree has trivial tail. 

To give a bit of history, for several years we conjectured in seminar talks 
that the RTP with marginal $\nu$ is non-endogenous. 
Because the simulation results suggested one of the condition 
equivalent to endogeny from \cite[Theorem 11]{AlBa05} fails
for the solution $\nu$ of the RDE (\ref{Equ:FPP-RDE}). 
In recent days for some time we thought we can
prove the opposite, but it turned out that our argument had some  
flaw in it. Fresh simulations confirm our earlier belief that the RTP
with marginal $\nu$ is non-endogenous. Till date to best of our knowledge
a rigorous proof is yet to be found. 

It is interesting to note that if the RTP with marginal $\nu$ is 
non-endogenous then the frozen percolation process would have a kind of
``spatial chaos'' property, that the behavior near the root would be affected 
by the behavior at infinity. On the other hand in light of the 
Theorem \ref{Thm:FPP-Tail-Trivial}, we note that possible influence of 
infinity at the root is not coming from the tail of the process. Such
examples are rare, our Example \ref{Ex:MC} is one such. But so far we
do not know a non-trivial example of this kind. Of course if 
non-endogeny for frozen percolation is proved, then that together with
Theorem \ref{Thm:FPP-Tail-Trivial} will provide one such.

\subsection{Outline of the Rest of the Paper} 
\label{SubSec:Outline}
The rest of the article is divided as follows. In the following section we
provide some basic connection between the root variable 
$X_{\emptyset}$ of an RTP with the tail $\sigma$-algebra $\HH$, 
and also give proofs of 
Propositions \ref{Prop:Endo} and \ref{Prop:Ex-MC}. 
In Section \ref{Sec:Proof-Main} we give a 
proof of the equivalence theorem and Section \ref{Sec:Proof-FPP-Tail}
contains the proof of the Theorem \ref{Thm:FPP-Tail-Trivial}. We conclude with
Section \ref{Sec:Remarks} which contains some further discussion.

\section{Connection between Root and Tail of a RTP}
\label{Sec:Root-Tail}
Because of the recursive structure one would expect that the tail 
$\sigma$-algebra $\HH$
is trivial, if and only if the root variable $X_{\emptyset}$ 
is independent of it. The following lemma preciously states that. 
\begin{Lemma}
\label{Lem:Tail-equiv}
$X_{\emptyset}$ is independent of $\HH$, if and only if $\HH$ is
trivial. 
\end{Lemma}

\proof If the tail $\HH$ is trivial then naturally 
$X_{\emptyset}$ is independent of it. 
For proving the converse we will need the following 
standard measure theoretic fact whose proof is a straightforward 
application of  
Dynkin's $\pi$-$\lambda$ theorem \cite{Bill95}, so we omit it here.  
\begin{Lemma}
\label{Lem:Tech-inde}
Suppose $\left(\Omega, \II, \bP \right)$ be a probability space and
let $\FF^*, \GG^*$ and $\HH^*$ be three sub-$\sigma$-algebras such that 
$\FF^*$ is independent of $\HH^*$;
$\GG^*$ is independent of $\HH^*$; and
$\FF^*$ and $\GG^*$ are independent given $\HH^*$. 
Then $\sigma\left(\FF^* \cup \GG^* \right)$ is independent of $\HH^*$. 
\end{Lemma}
To complete the proof of the Lemma \ref{Lem:Tail-equiv} we denote
$\FF_n^0 := \sigma\left( X_{\bi}, \vert \bi \vert = n \right)$ and
$\FF_n := \sigma\left( X_{\bi}, \vert \bi \vert \leq n \right)$.
From assumption $X_{\bi}$ is independent of $\HH$ for all $\bi \in \VV$.
Fix $n \geq 1$ and let $\bi \neq \bi'$ be two vertices at generation $n$.
From the definition
of RTP $X_{\bi}$ and $X_{\bi'}$ are independent,
moreover they are independent given $\HH_{n+k}$
for any $k \geq 1$. Letting $k \rightarrow \infty$ we conclude that 
$X_{\bi}$ and $X_{\bi'}$ are independent given $\HH$. Thus by
Lemma \ref{Lem:Tech-inde} we get that $\left( X_{\bi}, X_{\bi'} \right)$
is independent of $\HH$, and hence by induction $\FF_n^0$ is independent
of $\HH$. 

Now let 
$\GG_n = \sigma\left(\left\{ \left(\xi_{\bi}, N_{\bi}\right) \,\Big\vert\,
\vert \bi \vert \leq n \,\right\} \right)$, then $\GG_n$
is independent of $\HH$ from definition. Further 
$\GG_n$ is independent of $\FF_{n+1}^0$ 
given $\HH_{n+k}$ for any $k \geq 1$. 
Once again letting $k \rightarrow \infty$ 
we conclude that $\GG_n$ and $\FF_{n+1}^0$ are independent given
$\HH$. So again using Lemma \ref{Lem:Tech-inde} it follows that 
$\sigma\left( \GG_n \cup \FF_{n+1}^0 \right)$ is independent of
$\HH$. But 
$\FF_n \subseteq \sigma\left( \GG_n \cup \FF_{n+1}^0 \right)$
so $\FF_n$ is independent of $\HH$. 
But $\FF_n \uparrow \HH_0$ 
and hence $\HH$ is independent of $\HH_0 \supseteq \HH$. This proves
that $\HH$ is trivial. $\qed$

\subsection{Proof of Proposition \ref{Prop:Endo}}
\label{SubSec:Proof-of-Prop-Endo}
Let $\GG_n := \sigma \left( \left(\xi_{\bi}, N_{\bi}\right), 
\vert \bi \vert \leq n \right)$. 
From definition we have $\HH_n \downarrow \HH$ and 
$\GG_n \uparrow \GG$. Also for each $n \geq 0$, $\GG_n$ is
independent of $\HH_{n+1}$. So clearly $\GG$ is independent of $\HH$. 
Hence if the RTP is endogenous then 
$X_{\emptyset}$ is $\GG$-measurable, so it is independent of $\HH$. The
rest follows from the Lemma \ref{Lem:Tail-equiv}. \qed

\subsection{Proof of Proposition \ref{Prop:Ex-MC}}
\label{SubSec:Proof-of-Prop-Ex-MC}
There are several ways one can prove Proposition \ref{Prop:Ex-MC}, perhaps the
simplest is to apply the equivalence theorem (Theorem \ref{Thm:Bi-equiv-2}). 
This will also illustrate an easy application of the equivalence theorem. 
A non-trivial application is given in Section \ref{Sec:Bi-Uni-2nd-FPP} and
\ref{Sec:Proof-FPP-Tail}. 

\proof We will show that the bivariate uniqueness of the second kind holds
for the unique solution $\mbox{Bernoulli}\left(\sfrac{1}{2}\right)$ of the
RDE (\ref{Equ:RDE-Ex-MC}). So by part (b) of the equivalence theorem
(Theorem \ref{Thm:Bi-equiv-2}) the tail-triviality will follow (note
that in this case the continuity condition trivially holds). 

Let $\left(X,Y\right)$ be $S^2$-valued random pair with some
distribution such that the marginals are both Bernoulli$(1/2)$.
Let 
\[
\theta = \bP\left(X=1,Y=1\right) = \bP\left(X=0,Y=0\right).
\]
Suppose further that the distribution of 
$\left(X,Y\right)$ satisfies the following bivariate RDE
\[
\left(\begin{array}{c}
                      X \\
                      Y \end{array} \right) \ed
\left(\begin{array}{c}
                      X_1 + \xi \\
                      Y_1 + \eta \end{array} \right)
\,\,\,\, \left( \mbox{\ mod\ } 2 \,\right),
\]
where $\left(X_1, Y_1\right)$ is a copy of $\left(X,Y\right)$ and
independent of $\left(\xi, \eta\right)$ which are i.i.d.
Bernoulli$(q)$. So we get the following equation for $\theta$
\begin{equation}
\theta =  q^2 \theta + (1-q)^2 \theta + 2 q (1-q) (1/2 - \theta).
\label{Equ:Theta}
\end{equation}
The only solution of (\ref{Equ:Theta})
is $\theta=1/4$, thus $X$ and $Y$ must be independent, proving
the bivariate uniqueness of the second kind. \qed

\section{Proof of the Equivalence Theorem}
\label{Sec:Proof-Main}
(a) Let $\lambda$ be a fixed point of $\Tprod$ with
marginals $\mu$. 
Consider two independent and identical copies innovation processes
given by  $\left(\left( \xi_{\bi}, N_{\bi} \right), \bi \in \VV \right)$ and
$\left(\left( \eta_{\bi}, M_{\bi} \right), \bi \in \VV \right)$.  
Using Kolmogorov's 
consistency theorem \cite{Bill95},
we can then construct a bivariate RTP
$\left(\left(X^{(1)}_\bi,X^{(2)}_\bi\right), \bi \in \VV\right)$
with $\lambda = {\rm dist}(X^{(1)}_\emptyset,X^{(2)}_\emptyset)$. 
We note that this construction is no different than what one does
to obtain an univariate RTP as in (\ref{Equ:RTP}), and the bivariate RTP
has the similar properties as well. 
Notice that $\left(X^{(1)}_{\bi}\right)_{\bi \in \VV}$ and 
$\left(X^{(2)}_{\bi}\right)_{\bi \in \VV}$
are two (\emph{univariate}) RTPs with marginal $\mu$. So from assumption
both has trivial tails. 

We define the following $\sigma$-algebras
\begin{eqnarray}
\HH_n^{(1)} & := & 
\sigma\left( \left\{ X_{\bi}^{(1)} \, \Big\vert \, \vert \bi \vert \geq n 
\right\} \right) ; \\
\HH_n^{(2)} & := & 
\sigma\left( \left\{ X_{\bi}^{(2)} \, \Big\vert \, \vert \bi \vert \geq n 
\right\} \right) ; \\
\HH_n^{(*)} & := & 
\sigma\left( \left\{ \left(X_{\bi}^{(1)}, X_{\bi}^{(2)}\right) 
\, \Big\vert \, \vert \bi \vert \geq n  \right\} \right),
\end{eqnarray}
and we also define
\begin{eqnarray}
\mbox{Tail of\ } \left( X_{\bi}^{(1)} \right)_{\bi \in \VV} := \HH^{(1)}
& = & \mathop{\cap}\limits_{n \geq 0} \HH_n^{(1)}; \\
\mbox{Tail of\ } \left( X_{\bi}^{(2)} \right)_{\bi \in \VV} := \HH^{(2)}
& = & \mathop{\cap}\limits_{n \geq 0} \HH_n^{(2)}; \\
\mbox{Tail of\ } \left( X_{\bi}^{(1)}, X_{\bi}^{(2)} \right)_{\bi \in \VV} 
:= \HH^{(*)}
& = & \mathop{\cap}\limits_{n \geq 0} \HH_n^{(*)}.
\end{eqnarray}

Let $f$ and $g$ be two bounded measurable functions. Fix $n \geq 0$, 
\begin{eqnarray}
 &  & 
\bE\left[ f\left(X_{\emptyset}^{(1)}\right) g\left(X_{\emptyset}^{(2)}\right)
\, \Big\vert \, \HH_n^{(*)} \right] \nonumber \\
 & = & \bE\left[ f\left(X_{\emptyset}^{(1)}\right) \,\Big\vert\, \HH_n^{(*)} \right]
       \times 
       \bE\left[ g\left(X_{\emptyset}^{(2)}\right) \,\Big\vert\, \HH_n^{(*)} \right]
       \nonumber \\
 & = & \bE\left[ f\left(X_{\emptyset}^{(1)}\right) \,\Big\vert\, \HH_n^{(1)} \right]
       \times
       \bE\left[ g\left(X_{\emptyset}^{(2)}\right) \,\Big\vert\, \HH_n^{(2)} \right],
\end{eqnarray}
where the first equality follows from the recursive construction and 
because the two innovation processes are independent. 
Taking limit as $n \rightarrow \infty$ and using the martingale convergence
theorem we get
\begin{equation}
\bE\left[ f\left(X_{\emptyset}^{(1)}\right) g\left(X_{\emptyset}^{(2)}\right)
\,\Big\vert\, \HH^{(*)} \right]
= \bE\left[ f\left(X_{\emptyset}^{(1)}\right) \,\Big\vert\, \HH^{(1)} \right]
  \times 
  \bE\left[ g\left(X_{\emptyset}^{(2)}\right) \,\Big\vert\, \HH^{(2)} \right].
\end{equation}
Because both $\HH^{(1)}$ and $\HH^{(2)}$ are trivial,
so taking a further expectation we conclude that
\begin{equation}
\bE\left[ f\left( X_{\emptyset}^{(1)}\right) g\left( X_{\emptyset}^{(2)}\right)
\right] = \bE\left[ f\left(X_{\emptyset}^{(1)}\right) \right] \times
          \bE\left[ g\left(X_{\emptyset}^{(2)}\right) \right].
\end{equation}
So $X_{\emptyset}^{(1)}$ and $X_{\emptyset}^{(2)}$ are independent, that 
is, $\lambda = \mu \otimes \mu$, which implies that the bivariate uniqueness 
property of the second kind holds. 

(b) Let $(X_{\bi})_{\bi \in \VV}$ be the invariant RTP with marginal $\mu$.
$\HH_n$ and $\HH$ be as defined in (\ref{Equ:RTP-Tail-n}) and
(\ref{Equ:RTP-Tail}) respectively. 
Observe that $\HH_n \downarrow \HH$. 
Now fix $\Lambda : S \rightarrow \Rbold$ a bounded continuous
function. So by reverse martingale convergence 
\begin{equation} 
\bE \left[ \Lambda ( X_{\emptyset} ) \Big\vert \HH_n \right] 
\mathop{\longrightarrow}\limits^{\mbox{a.s.}}_{{\mathcal L}_2}
\bE \left[ \Lambda ( X_{\emptyset} ) \Big \vert \HH \right].
\label{martingale-limit-2}
\end{equation} 
Let $\left(\eta_{\bi}, M_{\bi}\right)_{\bi \in \VV}$ be independent
innovations which are independent of 
$\left(X_\bi\right)_{\bi \in \VV}$ and
$\left(\xi_{\bi}, N_{\bi}\right)_{\bi \in \VV}$.
For $n \geq 1$, define $Y_{\bi}^n := X_{\bi}$ if
$\vert \bi \vert=n$, and then recursively define $Y_{\bi}^n$
for $ \vert \bi \vert < n$ using RTP construction (\ref{Equ:RTP}),
but replacing $\xi_{\bi}$ by 
$\eta_{\bi}$ and $N_{\bi}$ by $M_{\bi}$
to get an invariant RTP $(Y^n_\bi)$ of depth $n$. 
Observe that $X_{\emptyset} \ed Y_{\emptyset}^n$. Further
given $\HH_n$, the variables $X_{\emptyset}$ and $Y_{\emptyset}^n$ 
are conditionally independent and identically distributed. 
Now let 
\begin{equation}
\bar{\sigma}_n^2 ( \Lambda ) :=
\Big\| \bE \left[ \Lambda ( X_{\emptyset} ) \Big\vert \HH_n \right] 
- \bE \left[\Lambda ( X_{\emptyset} ) \right] \Big\|_2^2 \,\, .
\label{sigma-def-2}
\end{equation}
We calculate
\begin{eqnarray}
\bar{\sigma}_n^2 ( \Lambda ) & = & \bE \left[ \left( 
                             \bE \left[ \Lambda ( X_{\emptyset} )
                               \Big\vert \HH_n \right]
                             - \bE\left[\Lambda(X_{\emptyset})\right]  
                             \right)^2 \right] \nonumber \\
                       & = & \var \left( \bE \left[ 
                             \Lambda ( X_{\emptyset} ) \Big\vert \HH_n 
                             \right] \right) \nonumber \\
                       & = & \var\left(\Lambda(X_{\emptyset})\right) -
                             \bE\left[
                             \var\left(\Lambda(X_{\emptyset})
                             \Big\vert \HH_n \right) \right] \nonumber \\
                       & = & \var\left(\Lambda(X_{\emptyset})\right) -
                             \sfrac{1}{2} 
                             \bE \left[ \left( \Lambda ( X_{\emptyset} ) 
                             - \Lambda  (Y_{\emptyset}^n )  \right)^2
                             \right]. \label{sigma-working-def-2} 
\end{eqnarray}
The last equality uses the conditional form of the fact that for any random 
variable $U$, one has $\var(U) = \sfrac{1}{2} \bE\left[ (U_1 - U_2)^2 \right]$,
where $U_1, U_2$ are i.i.d copies of $U$. 

Now suppose we show that 
\begin{equation}
 ( X_{\emptyset}, Y_{\emptyset}^n ) \cd
(X^{\star}, Y^{\star})
\label{toshow-2}
\end{equation}
for some limit $\left( X^{\star}, Y^{\star} \right)$.
From the construction, 
\[ 
\left[ \begin{array}{c}
                       X_{\emptyset} \\
                       Y_{\emptyset}^{n+1} \end{array} \right]
\ed
\left(\Tprod\right) \left( 
\left[ \begin{array}{c}
                       X_{\emptyset} \\
                       Y_{\emptyset}^n \end{array} \right]
\right),
\]
and then the weak continuity assumption on $\Tprod$ implies
\[ 
\left[ \begin{array}{c}
                       X^{\star} \\
                       Y^{\star} \end{array} \right]
\ed
\left(\Tprod\right) \left(
\left[ \begin{array}{c}
                       X^{\star} \\
                       Y^{\star} \end{array} \right]
\right).
\]
Also by construction we have 
$X_{\emptyset} \ed Y_{\emptyset}^n \ed \mu$ for all  $ n \geq 1$, and hence
$X^{\star} \ed Y^{\star} \ed \mu$. Now since we assume that the 
bivariate uniqueness property of the second kind holds, so 
$X^{\star}$ and $Y^{\star}$ must be independent.
Since $\Lambda$ is a bounded
continuous function, (\ref{toshow-2}) implies
\begin{equation}
\label{var-cov-2}
\bE\left[\left(\Lambda(X_{\emptyset}) - \Lambda(
    Y_{\emptyset}^n)\right)^2\right]  \rightarrow 
\bE\left[\left(\Lambda(X^{\star}) - \Lambda(
    Y^{\star})\right)^2\right] = 2
    \var\left(\Lambda(X_{\emptyset})\right)
\end{equation}
and so using (\ref{sigma-working-def-2}) we see that
$\bar{\sigma}_n^2 ( \Lambda ) \longrightarrow 0$. 
Hence from (\ref{sigma-def-2}) and (\ref{martingale-limit-2}) we conclude that
$\Lambda ( X_{\emptyset} )$ is independent of $\HH$.  This is true for every 
bounded continuous $\Lambda$, proving that
$X_{\emptyset}$ is independent of $\HH$, so from Lemma
\ref{Lem:Tail-equiv} it follows that $\HH$ is trivial. 

Now all remains is to show that limit (\ref{toshow-2}) exists.
Fix 
$f : S \rightarrow \Rbold$ and $h : S \rightarrow \Rbold$, two bounded
continuous functions. Again by reverse martingale convergence 
\[  
\bE \left[ f ( X_{\emptyset} ) \Big\vert \HH_n \right] 
\mathop{\longrightarrow}\limits^{\mbox{a.s.}}_{{\mathcal L}_1}
\bE \left[ f ( X_{\emptyset} ) \Big\vert \HH \right], 
\]
and similarly for $h$.
So
\begin{eqnarray*}
\bE \left[ f ( X_{\emptyset} ) h ( Y_{\emptyset}^n ) \right] & = &
\bE \left[ \bE \left[ f ( X_{\emptyset} ) h ( Y_{\emptyset}^n ) \Big\vert
\HH_n \right] \right] \\
 & = & \bE \left[ \bE \left[ f(X_{\emptyset}) \Big\vert \HH_n \right]
                  \bE \left[ h(X_{\emptyset}) \Big\vert \HH_n \right] \right], 
\end{eqnarray*}
the last equality because of conditional on $\HH_n$ 
$X_{\emptyset}$ and $Y_{\emptyset}^n$ are independent and identically 
distributed. 
Letting $n \to \infty$ we get 
\begin{equation}
\bE \left[ f ( X_{\emptyset} ) h ( Y_{\emptyset}^n ) \right]
\longrightarrow
\bE \left[ \bE \left[ f(X_{\emptyset}) \Big\vert \GG \right]
                  \bE \left[ h(X_{\emptyset}) \Big\vert \GG \right] \right].
\label{weak-conv-2}
\end{equation}
Moreover note that $X_{\emptyset} \ed Y_{\emptyset}^n \ed \mu$
and so the
sequence of bivariate distributions 
$ ( X_{\emptyset}, Y_{\emptyset}^n ) $ is
tight.   Tightness, together with convergence (\ref{weak-conv-2}) for all
bounded continuous $f$ and $h$, implies weak convergence of 
$( X_{\emptyset}, Y_{\emptyset}^n ) $ .

(c) First assume that 
$\left(\Tprod\right)^n \left( \mudiag \right) \cd \muprod$, then with
the same construction as done in part (b) we get that 
\[
\left( X_{\emptyset}, Y_{\emptyset}^n \right) \cd 
\left( X^{\star}, Y^{\star} \right),
\]
where $X^{\star}$ and $Y^{\star}$ are independent copies of $X_{\emptyset}$.
Further recall that $\Lambda$ is bounded continuous, thus using
(\ref{sigma-working-def-2}), (\ref{sigma-def-2}) and
(\ref{martingale-limit-2}) we conclude that 
$\Lambda(X_{\emptyset})$ is independent of $\HH$. Since it is true for any
bounded continuous function $\Lambda$, thus $X_{\emptyset}$ is
independent of $\HH$. Thus again by Lemma \ref{Lem:Tail-equiv} the RTP
has trivial tail. 

Conversely, suppose that the invariant RTP 
with marginal $\mu$ has trivial tail.
Let $\Lambda_1$ and $\Lambda_2$ be two bounded continuous
functions. Note that the variables 
$\left( X_{\emptyset}, Y_{\emptyset}^n \right)$, as defined in part (b)
has joint distribution $\left(\Tprod\right)^n \left( \mudiag \right)$.
Further, given $\HH_n$, they are conditionally independent and have
same conditional law as of $X_{\emptyset}$ given $\HH_n$. So 
\begin{eqnarray*}
\bE \left[ \Lambda_1(X_{\emptyset}) \, \Lambda_2(Y_{\emptyset}^n) \right] 
& = & 
\bE \left[ \bE\left[ \Lambda_1(X_{\emptyset}) \Big\vert \HH_n \right] \,
           \bE\left[ \Lambda_2(X_{\emptyset}) \Big\vert \HH_n \right]
           \right] \\
& \rightarrow &
\bE \left[ \bE\left[ \Lambda_1(X_{\emptyset}) \Big\vert \HH \right] \,
           \bE\left[ \Lambda_2(X_{\emptyset}) \Big\vert \HH \right]
           \right] \\
& = &
\bE \left[ \Lambda_1(X_{\emptyset}) \right]
\bE \left[ \Lambda_2(X_{\emptyset}) \right].
\end{eqnarray*}
The convergence is by reverse martingale convergence, and the last equality is
by tail triviality and Lemma \ref{Lem:Tail-equiv}. So from definition we get 
\[
\left(\Tprod\right)^n \left( \mudiag \right) 
\ed
\left( X_{\emptyset}, Y_{\emptyset}^n \right)
\cd 
\muprod .
\] 
$\qed$

\section{Bivariate Uniqueness Property of the Second Kind for the Frozen 
Percolation RDE}
\label{Sec:Bi-Uni-2nd-FPP}
In this section we prove the bivariate uniqueness property of the second
kind for the frozen percolation RDE (\ref{Equ:FPP-RDE}). 

\begin{Theorem}
\label{Thm:FPP-biuni}
Consider the following bivariate RDE,
\begin{equation}
\left( \begin{array}{c} X \\ Y \end{array} \right)
\ed
\left( \begin{array}{c} 
                    \Phi \left( X_1 \wedge X_2 ; U \right) \\
                    \Phi \left( Y_1 \wedge Y_2 ; V\right)
       \end{array} \right),
\label{Equ:FPP-biRDE}
\end{equation}
where $\left(X_j, Y_j\right)_{j = 1,2}$ are i.i.d with same joint law as
$\left(X, Y\right)$ and have same marginal distribution $\nu$ given by
\begin{equation}
\nu\left(dx\right) = \sfrac{dx}{2 x^2}, \,\,\,\, \sfrac{1}{2} < x < 1 ; 
\,\,\,\, \nu\left( \left\{ \infty \right\} \right) = \sfrac{1}{2} \, , 
\label{Equ:Defin-nu}
\end{equation}
and are independent of $\left(U, V\right)$ which are 
i.i.d. with $\mbox{Uniform}[0,1]$ distribution;
and $\Phi$ is given by (\ref{Equ:FPP-Phi-defi}). 
Then the unique solution of this bivariate RDE (\ref{Equ:FPP-biRDE}) is
the product measure $\nu \otimes \nu$. 
\end{Theorem}

\proof Since $\nu$ is a solution of the RDE (\ref{Equ:FPP-RDE}), so by
Lemma \ref{Lem:Bi-Obvious-Sol}(a), the product measure $\nu \otimes \nu$ is a
solution of the bivariate RDE (\ref{Equ:FPP-biRDE}). We will show
it is the unique solution. Suppose $\left(X, Y\right)$ is a solution of
(\ref{Equ:FPP-biRDE}), and 
let $F(x,y) := \bP\left( X \leq x, Y \leq y \right)$, for $x, y \in [0,1]$
be the joint distribution function. Notice that if
$(x,y) \in [0,1]^2 \setminus D$ where $D := \left[\sfrac{1}{2}, 1\right]^2$
then $F(x,y)=0$. Now from equation (\ref{Equ:FPP-biRDE}) if
$x, y \in \left[\sfrac{1}{2},1\right]$ then
\begin{eqnarray}
       &   & F(x,y) \nonumber \\
       & = & \bP\left( \Phi(X_1 \wedge X_2 ; U) \leq x, \,
                       \Phi(Y_1 \wedge Y_2 ; V) \leq y \right)
                       \nonumber \\
       & = & \bP\left( U < X_1 \wedge X_2 \leq x, \,
                       V < Y_1 \wedge Y_2 \leq y \right) 
             \nonumber \\
       & = & \bE\left[ \left( 1_{(X_1 \wedge X_2 > U)} - 
                              1_{(X_1 \wedge X_2 > x)} \right)
                       \left( 1_{(Y_1 \wedge Y_2 > V)} -
                              1_{(Y_1 \wedge Y_2 > y)} \right) 
                       \, 1_{(U < x)} \, 1_{(V < y)} \right] 
                       \nonumber \\
       & = & \int_0^x \int_0^y  \!
             \left( G^2(x,y) - G^2(x,v) - G^2(u,y) + G^2(u,v) \right) 
             \, dv \, du \label{Equ:Frozen-Int-1}
\end{eqnarray}
where $G(x,y) := \bP\left( X > x, Y> y \right)$, which can be written
as 
\begin{eqnarray}
G(x,y) & = & F(x,y) - \bP\left(X \leq x\right) 
                    - \bP\left( Y \leq y \right) + 1 \nonumber \\
       & = & F(x,y) + \sfrac{1}{2 x} + \sfrac{1}{2 y} - 1 \, .
\label{Equ:Frozen-FG}
\end{eqnarray}
Further notice that $G(x,y)=1$ if $x, y \leq \sfrac{1}{2}$; 
$G(x,y) = \frac{1}{2x}$ if $x \in \left[\sfrac{1}{2},1\right]$ 
and $y \leq \sfrac{1}{2}$; and finally
$G(x,y) = \frac{1}{2y}$ if $y \in \left[\sfrac{1}{2},1\right]$ 
and $x \leq \sfrac{1}{2}$. So (\ref{Equ:Frozen-Int-1}) can be written
as
\begin{eqnarray} 
F(x,y) & = & x y \, G^2(x,y) - \sfrac{1}{4 x} - \sfrac{1}{4 y} + \sfrac{3}{4} 
             \nonumber \\
       &   & - x \, \int_{\sfrac{1}{2}}^y \! G^2(x,v) \, dv
             - y \, \int_{\sfrac{1}{2}}^x \! G^2(u,y) \, du
             + \int_{\sfrac{1}{2}}^x \int_{\sfrac{1}{2}}^y \!
               G^2(u,v) \, dv \, du  \, , \nonumber \\
       &   & \,  \label{Equ:Frozen-Int-1.5}
\end{eqnarray}
when $x, y \in \left[ \sfrac{1}{2}, 1 \right]$.
We know that $G_0\left(x, y\right) := \frac{1}{4 x y}$ on 
$\left[\sfrac{1}{2},1\right] \times \left[\sfrac{1}{2},1\right]$ is a solution 
of the equation (\ref{Equ:Frozen-Int-1.5}) which represent the 
$\nu \otimes \nu$ solution of the bivariate equation
(\ref{Equ:FPP-biRDE}). Let $F_0$ be the distribution function for this
solution. Note that for this solution 
$F_0(1,1) = G_0(1,1) = \sfrac{1}{4}$ is the mass at the point
$\left( \infty, \infty \right)$. 

Let $H(x,y) = 1 - G(x,y)/G_0(x,y)$, where $0 \leq x, y \leq 1$. 
Notice that $H \equiv 0$ on $[0,1]^2 \setminus D$. 
Moreover for $(x,y) \in D$, 
\begin{eqnarray*}
G(x,y) &  =   & \bP\left( X > x , Y > y \right) \\
       & \leq & \min\left( \bP\left( X > x \right), 
                           \bP\left( Y > y \right)\right) \\
       &  =   & \frac{1}{2 (x \vee y)} \\
       & \leq & \frac{1}{2 x y} = 2 \, G_0(x,y) \, ,
\end{eqnarray*}
where the last inequality follows because $\sfrac{1}{2} \leq x, y \leq 1$.
Thus $-1 \leq H(x,y) \leq 1$ for all $(x,y) \in D$. 
To prove the bivariate uniqueness all we need to show is $H \equiv 0$ on 
$D$. 

Recall that $G_0$ satisfy (\ref{Equ:Frozen-Int-1.5}), that is,
\begin{eqnarray*}
F_0(x,y) & = & x y \, G_0^2(x,y) 
               - \sfrac{1}{4 x} - \sfrac{1}{4 y} + \sfrac{3}{4} 
               \nonumber \\
         &   & - x \, \int_{\sfrac{1}{2}}^y \! G_0^2(x,v) \, dv
               - y \, \int_{\sfrac{1}{2}}^x \! G_0^2(u,y) \, du
               + \int_{\sfrac{1}{2}}^x \int_{\sfrac{1}{2}}^y \!
                 G_0^2(u,v) \, du \, dv  \, .
\end{eqnarray*}
Further by (\ref{Equ:Frozen-FG}) and definition of $H$ we have 
\[ F_0(x,y) - F(x,y) = G_0(x,y) - G(x,y) = G_0(x,y) H(x,y). \]
So using (\ref{Equ:Frozen-Int-1.5}) we get 
\begin{eqnarray} 
&   & G_0(x,y) H(x,y) \nonumber \\
& = & F_0(x,y) - F(x,y) \nonumber \\
& = & x y \left( G_0^2(x,y) - G^2(x,y) \right) 
      + \int_{\sfrac{1}{2}}^x \int_{\sfrac{1}{2}}^y \! 
      \left(G_0^2(u,v) - G^2(u,v) \right) \, du \, dv \nonumber \\
&   & - x \, \int_{\sfrac{1}{2}}^y \! \left(G_0^2(x,v) - G^2(x,v)\right) \, dv
      - y \, \int_{\sfrac{1}{2}}^x \! \left(G_0^2(u,y) - G^2(u,y)\right) \, du
      \, . \nonumber \\
&   & \, \label{Equ:Frozen-Int-2}
\end{eqnarray}
Observe that 
\begin{equation}
G_0^2 - G^2 = G_0^2 - G_0^2 \left(1 - H\right)^2 = 
G_0^2 \left( 2 H - H^2 \right) \, ,
\label{Equ:Frozen-Int-2.5}
\end{equation}
and also using $G_0(x,y) = \sfrac{1}{4xy}$ on $D$, we get 
\begin{eqnarray} 
 &   & G_0(x,y) H(x,y) - x y \, G_0^2(x,y) \left( 2 H(x,y) - H^2(x,y) \right)
       \nonumber \\
 & = & G_0(x,y) H(x,y) 
       \left( 1 - x y \, G_0(x,y) \left( 2 - H(x,y) \right) \right) 
       \nonumber \\
 & = & \sfrac{1}{4} G_0(x,y) H(x,y) \left( 2 + H(x,y) \right) \nonumber \\
 & = & H(x,y) \, \frac{\left(2 + H(x,y)\right)}{16 x y} \, .
       \label{Equ:Frozen-Int-2.75} 
\end{eqnarray} 
Thus using (\ref{Equ:Frozen-Int-2}), (\ref{Equ:Frozen-Int-2.5}) and
(\ref{Equ:Frozen-Int-2.75}) we conclude that for $(x,y) \in D$,
\begin{eqnarray}
 &   & H(x,y) \nonumber \\ 
 & = & \frac{16 x y}{2 + H(x,y)} \left[ 
       \int_{\sfrac{1}{2}}^x \int_{\sfrac{1}{2}}^y \! 
       G_0^2(u,v) H(u,v) \left(2 - H(u,v)\right) \, dv \, du \right.
       \nonumber \\ 
 &   & \,\,\,\, \left.
       - x \, \int_{\sfrac{1}{2}}^y \! 
       G_0^2(x,v) H(x,v) \left(2-H(x,v)\right) \, dv 
       - y \, \int_{\sfrac{1}{2}}^x \! 
       G_0^2(u,y) H(u,y) \left(2-H(u,y)\right) \, du 
       \right] \nonumber \\
 & = & \frac{x y}{2 + H(x,y)} \left[ 
       \int_{\sfrac{1}{2}}^x \int_{\sfrac{1}{2}}^y \! 
       \sfrac{1}{u^2 v^2} H(u,v) \left(2 - H(u,v)\right) \, dv \, du 
       \right. \nonumber \\
 &   & \,\,\,\, \left.
       - \sfrac{1}{x} \, \int_{\sfrac{1}{2}}^y \! 
       \sfrac{1}{v^2} H(x,v) \left(2-H(x,v)\right) \, dv 
       - \sfrac{1}{y} \, \int_{\sfrac{1}{2}}^x \! 
       \sfrac{1}{u^2} H(u,y) \left(2-H(u,y)\right) \, du 
       \right] \, . \label{Equ:Frozen-Int-Key} 
\end{eqnarray}

Fix $0 < \eps < \sfrac{1}{3}$ then there exists a partition 
$\sfrac{1}{2} = a_0 < a_1 < a_2 < \ldots < a_{k-1} < a_k = 1$ 
of $\left[\sfrac{1}{2}, 1\right]$ with equal lengths, such that
\begin{equation}
\int_{a_i}^{a_{i+1}} \int_{a_j}^{a_{j+1}} \! \sfrac{dv \, du}{u^2 v^2} 
+ 2 \int_{a_i}^{a_{i+1}} \! \sfrac{du}{u^2} 
+ 2 \int_{a_j}^{a_{j+1}} \! \sfrac{dv}{v^2}
< \eps  \,\,\,\, \forall \,\,\, 0 \leq i,j \leq k-1 \, , 
\label{Equ:Frozen-Esti-1}
\end{equation}
where $(x, y) \in D$. This we can do because the function
$s \mapsto \sfrac{1}{s^2}$ is a continuous decreasing function on
$\left[ \sfrac{1}{2}, 1 \right]$. 

Put $B_{i,j} := \left[a_i, a_{i+1} \right] \times \left[a_j, a_{j+1}\right]$ 
and let  
$\norm H \norm_{i,j} := \sup_{x, y \in B_{i j}} \vert H(x,y) \vert$, 
for $0 \leq i, j \leq k-1$. 
Start with $i=j=0$ and let $(x,y) \in B_{i,j}$, observe that from equation
(\ref{Equ:Frozen-Int-Key}) we have
\begin{eqnarray}
 &      & \vert H(x,y) \vert \nonumber \\ 
 &   =  & \left\vert \frac{x y}{2 + H(x,y)} \right\vert 
          \times \left\vert
          \int_{\sfrac{1}{2}}^x \int_{\sfrac{1}{2}}^y \! 
          \sfrac{1}{u^2 v^2} \, H(u,v) \left(2 - H(u,v)\right) \, dv \, du 
          \right. \nonumber \\
 &      & \,\,\,\, \left.
          - \sfrac{1}{x} \int_{\sfrac{1}{2}}^y \! 
          \sfrac{1}{v^2} \, H(x,v) \left(2-H(x,v)\right) \, dv 
          - \sfrac{1}{y} \int_{\sfrac{1}{2}}^x \! 
          \sfrac{1}{u^2} \, H(u,y) \left(2-H(u,y)\right) \, du 
          \right\vert \nonumber \\
 &   =  & \left\vert \frac{x y}{2 + H(x,y)} \right\vert 
          \times \left\vert 
          \int_{a_i}^x \int_{a_j}^y \! 
          \sfrac{1}{u^2 v^2} \, H(u,v) \left(2 - H(u,v)\right) \, dv \, du 
          \right. \nonumber \\
 &      & \,\,\,\, \left.
          - \sfrac{1}{x} \int_{a_j}^y \! 
          \sfrac{1}{v^2} \, H(x,v) \left(2-H(x,v)\right) \, dv 
          - \sfrac{1}{y} \int_{a_i}^x \! 
          \sfrac{1}{u^2} \, H(u,y) \left(2-H(u,y)\right) \, du 
          \right\vert \nonumber \\
 & \leq & \frac{3 x y \, \norm H \norm_{i,j}}{2 + H(x,y)} 
          \times \left[ 
          \int_{a_i}^x \int_{a_j}^y \! \sfrac{dv \, du}{u^2 v^2}
          + 2 \int_{a_j}^y \! \sfrac{dv}{v^2} 
          + 2 \int_{a_i}^x \! \sfrac{du}{u^2}
          \right] \nonumber \\
 & \leq & 3 \norm H \norm_{i,j} \times
          \left[
          \int_{a_i}^{a_{i+1}} \int_{a_j}^{a_{j+1}} \! 
          \sfrac{dv \, du}{u^2 v^2}
          + 2 \int_{a_i}^{a_{i+1}} \! \sfrac{du}{u^2} 
          + 2 \int_{a_j}^{a_{j+1}} \! \sfrac{dv}{v^2}
          \right] \label{Equ:Frozen-Esti-2}
\end{eqnarray}
where the last but one inequality follows because 
$(x,y) \in B_{i,j} \subseteq D = \left[ \sfrac{1}{2}, 1 \right]^2$, and
so $x, y \geq \sfrac{1}{2}$, and also because $1 \leq 2 - H \leq 3$ on $D$,
and the last inequality follows because $2 + H \geq 1$ on $D$.  
So from (\ref{Equ:Frozen-Esti-2}) we get 
\[
\norm H \norm_{i,j} \leq 3 \, \eps \, \norm H \norm_{i,j} \, .
\]
But we have chosen $\eps < \sfrac{1}{3}$, so we must have 
\[
H(x,y) = 0 \,\,\,\, \mbox{for all} \,\,\, (x,y) \in B_{i,j} \, .
\]
Now we do induction on two indices $i$ and $j$ in the following way. 
For every fixed $0 \leq l \leq k-1$ we start with $i=j=l$ and then
continue with $i \in \left\{ l, l+1, \ldots, k-1 \right\}$, and 
$j \in \left\{ l, l+1, \ldots, k-1 \right\}$, repeating the above
argument in each step.  This finally yields 
\[ H(x,y) = 0 \,\,\,\, \mbox{for all} \,\,\, (x,y) \in D \, ,\] 
which completes the proof. \qed

\section{Proof of Theorem \ref{Thm:FPP-Tail-Trivial}}
\label{Sec:Proof-FPP-Tail}

Now to prove the Theorem \ref{Thm:FPP-Tail-Trivial} we will use the
part (b) of our equivalence theorem (Theorem \ref{Thm:Bi-equiv-2}). 
The bivariate uniqueness of the second kind has been proved in 
Theorem \ref{Thm:FPP-biuni}, so it only remains to check the technical 
condition of Theorem \ref{Thm:Bi-equiv-2}(b). 

For that suppose 
$\nu_n^{(2)} \cd \nu^{(2)}$ where $\left\{ \nu_n^{(2)} \right\}_{n \geq 1}$
and $\nu^{(2)}$ are bivariate distributions on $I^2$ with marginals $\nu$. 
Let $F_n$ be the distribution function for $\nu_n^{(2)}$ and $F$ be that
for $\nu^{(2)}$. We define $G_n$ and $G$ in similar manner as done in
equation (\ref{Equ:Frozen-FG}). 
Following argument similar of derivation of the equation
(\ref{Equ:Frozen-Int-1}) we get that for 
$x, y \in \left[ \sfrac{1}{2}, 1 \right]$,
\[
T \otimes T \left(F_n\right) (x,y) 
       = \int_0^x \int_0^y \!
         \left( G_n^2(x,y) - G_n^2(x,v) - G_n^2(u,y) + G_n^2(u,v) \right) 
         \, dv \, du \,.
\]
The rest follows using the dominated convergence theorem. \qed

\section{Remarks and Complement}
\label{Sec:Remarks}

\subsection{Tail-Triviality and Long Range Independence}

Gamarnik et. al. \cite{Ga04} introduced the concept of 
\emph{long range independence} for some particular RDEs, similar concept
was also used in later works \cite{Ba05a, BaGa05}. Borrowing their idea
we define the long range independence property for an
invariant RTP as follows.
\begin{Definition}
Suppose $\left(X_{\bi}\right)_{\bi \in \VV}$ be an invariant RTP with
marginal $\mu$, then we will say that the \emph{long range independence
property} holds if
\begin{equation}
\lim_{d \rightarrow \infty}
\mathop{\sup_{x_{\bi} \in S}}\limits_{\vert \bi \vert = d}
\rho\left(
\mbox{dist}\left( X_{\emptyset} \,\Big\vert\, X_{\bi} = x_{\bi},
\, \vert \bi \vert = d \,\right) \,,\, \mu \right)
= 0 \, , 
\label{Equ:Long-Range-Inde}
\end{equation}
where $\rho$ is a metric for the weak convergence topology on
$\PP\left(S\right)$.
\end{Definition}

\begin{Proposition}
\label{Lem:Lang-to-Tail}
Suppose $\left(X_{\bi}\right)_{\bi \in \VV}$ is an invariant RTP with
marginal $\mu$ which has long range independence
property as defined above, then it must have trivial tail. 
\end{Proposition}

\proof Let $\HH = \mathop{\cap}\limits_{n \geq 0} \HH_n$ be the tail of the RTP
$\left(X_{\bi}\right)_{\bi \in \VV}$ where $\HH_n$ is as defined
in (\ref{Equ:RTP-Tail-n}). Let $\Lambda : S \rightarrow \Rbold$ be a bounded
continuous function and consider the
conditional expectation 
$\bE\left[ \Lambda\left(X_{\emptyset}\right) \,\Big\vert\, \HH_n \right]$,
by martingale convergence theorem 
\[
\bE\left[ \Lambda\left(X_{\emptyset}\right) \,\Big\vert\, \HH_n \right]
\longrightarrow
\bE\left[ \Lambda\left(X_{\emptyset}\right) \,\Big\vert\, \HH \right]
\,\,\,\, \mbox{a.s.}
\]
On the other hand from the long range independence property it follows
\[
\bE\left[ \Lambda\left(X_{\emptyset}\right) \,\Big\vert\, \HH_n \right]
\longrightarrow
\bE\left[ \Lambda\left(X_{\emptyset}\right) \right]
\,\,\,\, \mbox{a.s.} \,,
\]
since $X_{\emptyset} \sim \mu$. Thus we get 
\[
\bE\left[ \Lambda\left(X_{\emptyset}\right) \,\Big\vert\, \HH \right]
=
\bE\left[ \Lambda\left(X_{\emptyset}\right) \right]
\,\,\,\, \mbox{a.s.} \,,
\]
which is true for every bounded continuous function $\Lambda$, hence 
we must have $X_{\emptyset}$ independent of $\HH$. So by 
Lemma \ref{Lem:Tail-equiv} we conclude that the tail of the RTP is
trivial. $\qed$

Now the converse is not necessarily true. To see this, we first note that
in order for a RTP to have the long range independence,
the underlying RDE 
need to satisfy certain properties. For example,
\begin{Lemma}
\label{Lem:Long-Unique}
Suppose an invariant RTP with marginal $\mu$ has long range independence 
property. If $T$ is the associated operator for the RDE with 
domain $\PP$ then for any $\mu' \in \PP$ we must have
\begin{equation}
T^n\left(\mu' \right) \cd \mu \,\,\,\, \mbox{as} \,\,\, 
n \rightarrow \infty \,.
\label{Equ:Domain-of-Attraction}
\end{equation}
\end{Lemma}
The proof of this lemma easily follows  
from equation (\ref{Equ:Long-Range-Inde}), the details are left for the
readers. But
from this lemma we see that if an invariant RTP with marginal $\mu$ has
long range independence property then the underlying RDE necessarily has
unique solution $\mu$. 
Now \cite{AlBa05} gives several examples of
RDEs which may have multiple solutions but some of which can be endogenous.
To give a specific example, we consider the Quicksort RDE, which is given
by
\begin{equation}
X \ed U X_1 + \left(1-U\right) X_2 + 
2 U \log U + 2 \left(1-U\right) \log \left(1 - U\right) + 1
\,\,\,\, \mbox{on} \,\,\, \Rbold, 
\label{Equ:Quicksort-RDE}
\end{equation}
It is known that this RDE has a two parameter family of solutions
\cite{FiJa00},  
and only those with finite first moment are endogenous 
(see \cite[Theorem 21]{AlBa05}).
So an invariant RTP with a marginal which is a solution of
(\ref{Equ:Quicksort-RDE}) and has finite first moment, 
will be endogenous and hence from Proposition \ref{Prop:Endo}
has trivial tail. 
But by Lemma \ref{Lem:Long-Unique} we conclude that this invariant RTP
can not have long range independence property because, the RDE 
(\ref{Equ:Quicksort-RDE}) has many solutions. 

Finally, even though
it is not quite related to tail-triviality, but we still note that 
the above example also shows that endogeny does not imply long range
independence property. Interesting enough the converse is not true either. 
It is in fact easy to show that the unique invariant RTP of the
Example \ref{Ex:MC}
discussed in Section \ref{Sec:Intro} 
has long range independence property, but it is not endogenous.
In light of
Lemma \ref{Lem:Long-Unique} one may conjecture that 
if a RDE has unique solution with full domain of attraction,
and the solution is endogenous, then it must have the long range
independence property, but this to best of our knowledge remains as an
open problem.

\subsection{Frozen Percolation on $r$-regular Trees}
\label{Subsec:FPP-d-regular}
Using exactly similar arguments as done in the case of infinite regular 
binary tree
one can construct an automorphism 
invariant version of frozen percolation process on a infinite
$r$-regular tree $\TTT_r$ in which each vertex has degree $r \geq 3$
(see \cite{Al00} for details). In this setting the RDE is given by
\begin{equation}
Y^r \ed \Phi^r \left( Y^r_1 \wedge Y^r_2 \wedge \cdots \wedge Y^r_{r-1} ; U
               \right) \,\,\, \mbox{on\ \ } 
                I^r := \left[ \sfrac{1}{r-1}, 1 \right] \cup \{ \infty \},
\label{Equ:FPP-r-RDE}
\end{equation}
where $\left(Y^r_j\right)_{1 \leq j \leq r-1}$ are i.i.d with same law
as $Y^r$ and are independent of $U \sim \mbox{Uniform}[0,1]$; and 
$\Phi^r : I^r \times [0,1] \rightarrow I^r$ 
is the function defined by equation (\ref{Equ:FPP-Phi-defi}). 
It is easy to check that the unique solution of this RDE with full support 
and having no atom in $\left[ \sfrac{1}{r-1}, 1 \right]$
is given by 
\begin{equation}
\nu^r\left( dy \right) = 
\frac{dy}{(r-2)(r-1)^{\sfrac{1}{r-2}} \,\, y^{\sfrac{r-1}{r-2}}},
\,\,\, \sfrac{1}{r-1} < y < 1, \,\,\,\,
\nu^r(\{ \infty \}) = \frac{1}{(r-1)^{\sfrac{1}{r-2}}}.
\label{Equ:FPP-r-Sol}
\end{equation}
Naturally the case $r=3$ gives back the RDE (\ref{Equ:FPP-RDE}) and its
fundamental solution $\nu$. Interesting enough our argument to prove 
the bivariate uniqueness of the second kind for the frozen percolation 
RDE (\ref{Equ:FPP-RDE}) extend essentially unchanged in this setting 
(only the constants need to be changed). So the
invariant RTP associated with the RDE (\ref{Equ:FPP-r-RDE}) with 
marginal $\nu^r$ also has trivial tail. Once again the question of 
non-endogeny remains open.

\section*{Acknowledgment} 
An incomplete version of the main result of this article appears in 
author's doctoral dissertation, written under the guidance of 
Professor David J. Aldous, whom the author would like to thank for many 
discussion, and also for his continuous encouragement. 
The author would also like to thank Professor Don Aronson for his
insightful remarks on solving non-linear integral equations, and 
to Krishanu Maulik for carefully reading an earlier version of this paper. 
Thanks are also due to an anonymous referee for a through  
review of the paper. 

\bibliographystyle{plain}

\bibliography{RTP_Tail.bib}

\begin{thebibliography}{10}

\bibitem{Al92a}
David Aldous.
\newblock Asymptotics in the random assignment problem.
\newblock {\em Probab. Theory Related Fields}, 93(4):507--534, 1992.

\bibitem{AlSt04}
David Aldous and J.~Michael Steele.
\newblock The objective method: probabilistic combinatorial optimization and
  local weak convergence.
\newblock In {\em Probability on discrete structures}, volume 110 of {\em
  Encyclopaedia Math. Sci.}, pages 1--72. Springer, Berlin, 2004.

\bibitem{Al00}
David~J. Aldous.
\newblock The percolation process on a tree where infinite clusters are frozen.
\newblock {\em Math. Proc. Cambridge Philos. Soc.}, 128(3):465--477, 2000.

\bibitem{Al01}
David~J. Aldous.
\newblock The {$\zeta(2)$} limit in the random assignment problem.
\newblock {\em Random Structures Algorithms}, 18(4):381--418, 2001.

\bibitem{AlBa05}
David~J. Aldous and Antar Bandyopadhyay.
\newblock A survey of max-type recursive distributional equations.
\newblock {\em Ann. Appl. Probab.}, 15(2):1047--1110, 2005.

\bibitem{Ba02}
Antar Bandyopadhyay.
\newblock Bivariate {U}niqueness in the {L}ogistic {R}ecursive {D}istributional
  {E}quation.
\newblock Technical Report 629, Department of Statistics, UC Berkeley, 2002.

\bibitem{Ba05a}
Antar Bandyopadhyay.
\newblock Hard-{C}ore {M}odel on {R}andom {G}raphs.
\newblock (in preparation), 2005.

\bibitem{BaGa05}
Antar Bandyopadhyay and David Gamarnik.
\newblock Counting without sampling. {N}ew algorithms for enumeration problems
  using statistical physics.
\newblock To appear in ACM-SIAM Symposium on Discrete Algorithms 2006,
  (available at {\tt <http://www.arxiv.org/ps/math.PR/0510471>}), 2006.

\bibitem{Bill95}
Patrick Billingsley.
\newblock {\em Probability and measure}.
\newblock Wiley Series in Probability and Mathematical Statistics. John Wiley
  \& Sons Inc., New York, third edition, 1995.

\bibitem{FiJa00}
James~Allen Fill and Svante Janson.
\newblock A characterization of the set of fixed points of the {Q}uicksort
  transformation.
\newblock {\em Electron. Comm. Probab.}, 5:77--84 (electronic), 2000.

\bibitem{Ga04}
D.~Gamarnik, T.~Nowicki, and G.~Swirscsz.
\newblock Maximum {W}eight {I}ndependent {S}ets and {M}atchings in {S}parse
  {R}andom {G}raphs. {E}xact {R}esults using the {L}ocal {W}eak {C}onvergence
  {M}ethod.
\newblock To appear in \emph{Random Structures and Algorithms}, (available at
  {\tt <http://www.arxiv.org/pdf/math.PR/0309441>}), 2004.

\bibitem{RosRu01}
U.~R{\"o}sler and L.~R{\"u}schendorf.
\newblock The contraction method for recursive algorithms.
\newblock {\em Algorithmica}, 29(1-2):3--33, 2001.

\bibitem{Ros92}
Uwe R{\"o}sler.
\newblock A fixed point theorem for distributions.
\newblock {\em Stochastic Process. Appl.}, 42(2):195--214, 1992.

\end{thebibliography}


\end{document}